\documentclass[12pt]{amsart}
\usepackage{amsxtra}
\usepackage{amssymb}
\newcommand{\cleqn}{\setcounter{equation}{0}}
\newcommand{\clth}{\setcounter{theorem}{0}}
\newcommand {\sectionnew}[1]{\section{#1}\cleqn\clth}

\theoremstyle{plain}
\newtheorem{theorem}{Theorem}[section]
\newtheorem{lemma}[theorem]{Lemma}
\newtheorem{definition-theorem}[theorem]{Definition-Theorem}
\newtheorem{proposition}[theorem]{Proposition}
\newtheorem{corollary}[theorem]{Corollary}

\newtheorem{exampleproclaim}[theorem]{Example}
\newtheorem*{theoremA*}{Theorem A}
\newtheorem*{theoremB*}{Theorem B}
\newtheorem*{theoremC*}{Theorem C}
\theoremstyle{definition}

\newtheorem{remark}[theorem]{Remark}
\newtheorem*{remark*}{Remark}

\newtheorem*{problem*}{Problem}
\newcommand \NN {{\mathbb N}}

\newcommand \calR {{\mathcal{R}}}

\newcommand \Rtil {\widetilde{R}}

\newcommand \Ttil {\widetilde{T}}

\DeclareMathOperator\rann{r{.}ann}
\DeclareMathOperator\End{End}
\DeclareMathOperator\soc{soc}
\DeclareMathOperator\op{op}
\DeclareMathOperator\chr{char}
\DeclareMathOperator\rcfm{RCFM}
\DeclareMathOperator\fm{FM}

\def\top{T^{\op}}

\begin{document}
%%%%%%%%%%%%%%%%%%%%%%%%%%%%%%%%%%%%%%%%%%%%%%

\title{Extensions of locally matricial and locally semisimple algebras}
\author{K.R. Goodearl}
\address{
Department of Mathematics \\
University of California\\
Santa Barbara, CA 93106 \\
U.S.A.
}
\email{goodearl@math.ucsb.edu}

\begin{abstract}
Two extension problems are solved.  First, the class of locally matricial algebras over an arbitrary field is closed under extensions.  Second, the class of locally finite dimensional semisimple algebras over a fixed field is closed under extensions if and only if the base field is perfect.  Regardless of the base field, extensions of the latter type are always locally unit-regular.
\end{abstract}

\subjclass[2020]{}
\keywords{}

\maketitle

%%%%%%%%%%%%%%%%%%%%%%%%%%%%%
\sectionnew{Introduction}

An \emph{extension} of a ring (resp., algebra) $A$ by a ring (resp., algebra) $C$ is a ring (resp., algebra) $B$ which contains an ideal $J$ such that $J \cong A$ and $B/J \cong C$ as rings (resp., algebras).  This may be displayed as a short exact sequence
$$
0 \longrightarrow A \longrightarrow B \longrightarrow C \longrightarrow 0
$$
of ring (resp., algebra) homomorphisms.  The \emph{Extension Problem} for a class $\calR$ of rings or algebras is the problem of determining whether or not $\calR$ is closed under extensions.  Here we solve the Extension Problem for two classes of algebras defined by local properties, as follows.

Recall that an algebra $R$ over a field $K$ is \emph{matricial} provided $R \cong \prod_{i=1}^k M_{n_i}(K)$ for some finite list of positive integers $n_1,\dots,n_k$, and that $R$ is \emph{locally matricial} provided each finite subset of $R$ is contained in a matricial subalgebra of $R$.  Similarly, $R$ is a \emph{locally finite dimensional semisimple} $K$-algebra provided each finite subset of $R$ is contained in a finite dimensional semisimple subalgebra of $R$.  In particular, the term ``locally finite dimensional semisimple" means ``locally [finite dimensional semisimple]" as opposed to ``[locally finite dimensional] semisimple".  Finally, $R$ is a \emph{locally finite dimensional} $K$-algebra provided each finite subset of $R$ is contained in a finite dimensional subalgebra of $R$, or equivalently, provided each finitely generated subalgebra of $R$ is finite dimensional.  (This condition is often abbreviated to ``locally finite" in the literature.)  It has been known for more than six decades that the class of locally finite dimensional $K$-algebras is closed under extensions \cite[Proposition X.12.1]{Jac}.

\begin{theoremA*}  \label{thmA}
{\rm[Theorem \ref{lmxlm>lm}]}
Let $K$ be an arbitrary field.  Then any extension of a locally matricial $K$-algebra by a locally matricial $K$-algebra is locally matricial.
\end{theoremA*}

\begin{theoremB*}  \label{thmB}
{\rm[Theorem \ref{perf.lfdss.ext} and Example \ref{ex.charp}]}

{\rm(1)} Let $K$ be a perfect field.  Then any extension of a locally finite dimensional semisimple $K$-algebra by a locally finite dimensional semisimple $K$-algebra is locally finite dimensional semisimple.

{\rm(2)} Let $K$ be a non-perfect field.  Then there exists an extension of a locally finite dimensional semisimple $K$-algebra by a locally finite dimensional semisimple $K$-algebra which fails to be locally finite dimensional semisimple.
\end{theoremB*}

On the other hand, we exhibit classes of extensions of locally finite dimensional semisimple algebras over arbitrary fields which retain the property of being locally finite dimensional semisimple (Theorems \ref{ext.via.chekanu} and \ref{lfdssxLM>lfdss}).

We do not address the Extension Problem for locally semisimple rings or algebras here, except to mention a well known example of Menal and Moncasi \cite[Example 1]{MeMo}.  This is a unital regular $K$-algebra $R$ such that $\soc R$ is locally finite dimensional semisimple and $R/\soc R$ is a field while $R$ is directly infinite.  In particular, $R$ is not unit-regular, and so cannot be locally semisimple.

In the final section of the paper, we consider locally unit-regular rings and algebras.  Although the class of locally unit-regular $K$-algebras is not closed under extensions, as the Menal-Moncasi example shows, we prove that all extensions of locally finite dimensional semisimple by locally finite dimensional semisimple $K$-algebras are locally unit-regular.  More generally:

\begin{theoremC*}  \label{thmC}
{\rm[Theorem \ref{lfdssxlur>lur}]}
Let $K$ be an arbitrary field.  Then any extension of a locolly unit-regular $K$-algebra by a locally finite dimensional semisimple $K$-algebra is locally unit-regular.
\end{theoremC*}

\subsection{Conventions}
Fix a base field $K$ throughout.  \emph{Algebras} in this paper are assumed to be associative, but not necessarily unital, $K$-algebras.  A subalgebra which is unital as an algebra in its own right will be referred to as a \emph{subalgebra with identity}; there is no assumption that the ambient algebra is unital, nor that, if the ambient algebra is unital, its identity element agrees with the identity of the subalgebra.  \emph{Ideals} are assumed to be algebra ideals, i.e., (two-sided) ring ideals which are also $K$-subspaces.  \emph{Extensions} of one $K$-algebra by another $K$-algebra are assumed to be extensions within the category of $K$-algebras.

We often abbreviate \emph{finite dimensional} (meaning finite dimensional over $K$) to \emph{f.d}, and similarly abbreviate \emph{finite dimensional semisimple} to \emph{f.d.ss}.

%%%%%%%%%%%%%%%%%%%%%%%%%%%%%
\sectionnew{Regularity}  \label{sec.regular}

In effect, most of the work in this paper takes place in the universe of (von Neumann) regular rings, since (a) any locally matricial or locally finite dimensional semisimple algebra is regular, and (b) all extensions of regular rings by regular rings are regular (e.g., \cite[Lemma 1.3]{vnrr}).  In particular, any extension of a locally matricial or locally finite dimensional semisimple algebra by another such algebra is regular.

Any regular ring $R$ is locally unital: each finite subset of $R$ is contained in a corner $eRe$ for some idempotent $e \in R$ \cite[Lemma 2]{FaUt}.

Any ideal $J$ of a regular ring $R$ is regular, as is any corner $eRe$ with $e$ idempotent.  Since $eJe$ is an ideal of $eRe$ (whether or not $e \in J$), it follows that $eJe$ is regular.  Note that $eRe$ is an extension of $eJe$ by $(e+J)(R/J)(e+J)$.

We shall use the above facts about regular rings without further mention.

%%%%%%%%%%%%%%%%%%
\sectionnew{Locally matricial by locally matricial algebras}

In this section, Theorem A is proved.  In the process, we shall reduce to the case of an extension of a locally matricial algebra by a matricial algebra.  We first develop some intermediate results.

\begin{lemma}  \label{lm>corner}
Suppose $R$ is a locally matricial algebra, $J$ an ideal of $R$, and $e \in R$ an idempotent.  Then $eRe$ and $eJe$ are locally matricial algebras.
\end{lemma}

\begin{proof}
It suffices to deal with $eJe$, since the case of $eRe$ will follow by taking $J=R$.

Let $X$ be a finite subset of $eJe$.  Since $eJe$ is regular, there is an idempotent $f \in eJe$ such that $X \subseteq f(eJe)f = fJf$, and it suffices to show that $X$ is contained in a matricial subalgebra of $fJf$.
By hypothesis, there is a matricial subalgebra $T$ of $R$ that contains $X \cup \{f\}$.  Then $X \subseteq fTf \subseteq fRf = fJf$, so it is enough to show that $fTf$ is matricial.

Now $T = T_1 \oplus\cdots\oplus T_m$ where $T_i$ is an ideal of $T$ and $T_i \cong M_{n_i}(K)$ as algebras for some $n_i \in \NN$.  Write $f = f_1+\cdots+f_m$ for some idempotents $f_i \in T_i$.  As a right $T_i$-module, $f_iT_i$ is a finite direct sum of some $r_i$ copies of the unique simple right $T_i$-module $S_i$, and so $f_iTf_i \cong M_{r_i}(\End_{T_i}(S_i)) \cong M_{r_i}(K)$.  Thus $fTf \cong M_{r_1}(K) \oplus\cdots\oplus M_{r_m}(K)$, proving that $fTf$ is matricial, as required.
\end{proof}

\begin{lemma}  \label{R=T+J}
Let $R$ be a regular algebra, $J$ an ideal of $R$, and $\pi : R \rightarrow R/J$ the quotient map.  If $A$ is a finite dimensional subalgebra of $R/J$, then $R$ has a finite dimensional subalgebra $T$ such that $\pi$ maps $T$ isomorphically onto $A$.

In particular, if $R/J$ is finite dimensional, there exists a subalgebra $T$ of $R$ such that $R = T \oplus J$ as vector spaces.
\end{lemma}

\begin{proof}
Choose a finite dimensional subspace $A_0 \subseteq R$ such that $\pi(A_0) = A$, and an idempotent $e \in R$ such that $A_0 \subseteq eRe$.  Thus $A \subseteq \pi(e) (R/J) \pi(e)$, and it suffices to find a finite dimensional subalgebra $T$ of $eRe$ such that $T \cap eJe = 0$ and $T+eJe = e\pi^{-1}(A)e$.  After replacing $R$, $J$ and $A$ by $eRe$, $eJe$ and $e\pi^{-1}(A)e/eJe$, we may now assume that $R$ is unital.

Choose a basis $(a_1,\dots,a_n)$ for $A$.  There are scalars $\lambda^{(k)}_{ij} \in K$ such that $a_i a_j = \sum_k \lambda^{(k)}_{ij} a_k$ for all $i$, $j$.  Choose elements $x_i \in R$ such that $\pi(x_i) = a_i$, and set $v_{ij} := x_i x_j - \sum_k \lambda^{(k)}_{ij} x_k \in J$ for all $i$, $j$.

Now $V := \sum_{i,j} v_{ij} R$ contains all the $v_{ij}$ (by regularity), and $W := V + \sum_i x_iV$ is a finitely generated right ideal of $R$ contained in $J$, so $W = fR$ for some idempotent $f \in J$.  In particular, $v_{ij} = f v_{ij}$ for all $i$, $j$.  Given $i$, we have $x_iV \subseteq W$ and 
$$x_i x_j V = \bigl( \sum_k \lambda^{(k)}_{ij} x_k + v_{ij} \bigr) V \subseteq \sum_k x_kV + v_{ij}R \subseteq W \qquad \forall\; j,
$$
so $x_iW \subseteq W$.  Hence, $x_if = fx_if$.

Setting $g := 1-f$, we have $g v_{ij} = 0$ for all $i$, $j$ and $g x_i g = g x_i$ for all $i$.  Consequently, 
$$
gx_i gx_j = g x_ix_j = g \sum_k \lambda^{(k)}_{ij} x_k + g v_{ij} = \sum_k \lambda^{(k)}_{ij} gx_k \qquad \forall\; i,j.
$$
Thus $T := \sum_i K gx_i$ is a subalgebra of $R$.  Since $\pi(gx_i) = a_i$ for all $i$, the $gx_i$ are linearly independent and $\pi|_T$ is injective.  Therefore $\pi|_T$ gives an isomorphism of $T$ onto $A$.
\end{proof}

\begin{proposition}  \label{q.cent.T}
Let $R$ be a unital regular algebra, $J$ an ideal of $R$, and $T$ a matricial subalgebra of $R$ such that $R = T \oplus J$ as vector spaces.  Given a finite subset $Y \subseteq J$, there is an idempotent $e \in J$ such that $Y \subseteq eJe$ and $e$ centralizes $T$.
\end{proposition}

\begin{proof}
Set $p := 1_T$, and note that $1-p \in J$.  Hence, there is an idempotent $f \in J$ such that $Y \cup \{1-p\} \subseteq fJf$. Note that $f = pf+1-p$ where $pf$ is an idempotent in $pRp$.

Let $(e^{(s)}_{ij})$ for $s = 1,\dots,t$ and $i,j = 1,\dots,n_s$ be a full set of matrix units for $T$, hence also a $K$-basis.  The elements $p_s := \sum_{i=1}^{n_s} e^{(s)}_{ii}$ are central idempotents in $T$, with $\sum_{s=1}^t p_s = p$.

For each $r,s = 1,\dots,t$, $i = 1,\dots,n_s$, and $j = 1,\dots,n_r$, choose $u^{(rs)}_{ij} \in J$ such that
$$
(e^{(s)}_{1i} f e^{(r)}_{j1}) u^{(rs)}_{ij} (e^{(s)}_{1i} f e^{(r)}_{j1}) = e^{(s)}_{1i} f e^{(r)}_{j1}.
$$
After replacing $u^{(rs)}_{ij}$ by $e^{(r)}_{11} u^{(rs)}_{ij} e^{(s)}_{11}$, we may assume that $u^{(rs)}_{ij} \in e^{(r)}_{11} J e^{(s)}_{11}$. For each $s$, there is an idempotent $g_s \in e^{(s)}_{11} J e^{(s)}_{11}$ such that 
$$
e^{(s)}_{1i} f e^{(r)}_{j1} u^{(rs)}_{ij},\; u^{(sr)}_{ij} e^{(r)}_{1i} f e^{(s)}_{j1} \in g_s J g_s
$$
for all $r$ and all appropriate $i$, $j$.  Set $q_s := \sum_{i=1}^{n_s} e^{(s)}_{i1} g_s e^{(s)}_{1i}$, which is an idempotent in $p_sJp_s \subseteq pJp$ such that 
$$
e^{(s)}_{ij} q_s = e^{(s)}_{i1} g_s e^{(s)}_{1j} = q_s e^{(s)}_{ij} \qquad \forall\; i,j.
$$
It follows that $q_s e^{(s)}_{i1} g_s = e^{(s)}_{i1} q_s g_s = e^{(s)}_{i1} g_s e^{(s)}_{11} g_s = e^{(s)}_{i1} g_s$ for all $i$, and likewise $g_s e^{(s)}_{1j} q_s = g_s e^{(s)}_{1j}$ for all $j$.  Since
$$
e^{(s)}_{ii} f e^{(r)}_{j1} u^{(rs)}_{ij} = e^{(s)}_{i1} e^{(s)}_{1i} f e^{(r)}_{j1} u^{(rs)}_{ij} \in e^{(s)}_{i1} g_s J g_s \subseteq q_s J
$$
for all $r$, $i$, $j$, we obtain 
$$
e^{(s)}_{ii} f e^{(r)}_{j1} = e^{(s)}_{i1} e^{(s)}_{1i} f e^{(r)}_{j1} = e^{(s)}_{i1} (e^{(s)}_{1i} f e^{(r)}_{j1}) u^{(rs)}_{ij} (e^{(s)}_{1i} f e^{(r)}_{j1}) \in e^{(s)}_{ii} f e^{(r)}_{j1} u^{(rs)}_{ij} R \subseteq q_s J ,
$$
and thus $e^{(s)}_{ii} f e^{(r)}_{jj} \in q_s J$.  Similarly, $e^{(s)}_{ii} f e^{(r)}_{jj} \in J q_r$, whence $e^{(s)}_{ii} f e^{(r)}_{jj} \in q_s J q_r$.  Therefore
$$
p_s f p_r = \sum_{i,j} e^{(s)}_{ii} f e^{(r)}_{jj}  \in q_s J q_r \qquad \forall\; r,s.
$$

Since $q_1,\dots,q_t$ are pairwise orthogonal, $q := q_1 +\cdots+ q_t$ is an idempotent in $pJp$.  Now $p_s f p_r \in qJq$ for all $r$, $s$, and hence
$$
pf = pfp = \sum_{r,s=1}^t p_s f p_r \in qJq.
$$
Since $q_s$ commutes with $e^{(s)}_{ij}$ for all $i$, $j$, we see that $q$ commutes with  all $e^{(s)}_{ij}$, and therefore $q$ centralizes $T$.  

Finally, $e := q+1-p$ is an idempotent in $J$.  Since $1-p$ annihilates $T$ on either side, $e$ centralizes $T$.  Moreover, $f = pf+1-p \in eJe$ because $pf \in qJq$, and thus $Y \subseteq fJf \subseteq eJe$.
\end{proof}

\begin{theorem}  \label{lmxlm>lm}
Let $R$ be a $K$-algebra and $J$ an ideal of $R$.  If $J$ and $R/J$ are both locally matricial $K$-algebras, then $R$ is locally matricial.
\end{theorem}

\begin{proof}
Given a finite set $X \subseteq R$, we need to show that $X$ is contained in a matricial subalgebra of $R$.

There is an idempotent $e \in R$ such that $X \subseteq eRe$.  By Lemma \ref{lm>corner}, $eJe$ and $eRe/eJe$ are locally matricial, and it suffices to prove that $X$ is contained in a matricial subalgebra of $eRe$.  Thus, we may assume that $R$ is unital.  Now let $\pi : R \rightarrow R/J$ be the quotient map.  Since $R/J$ is locally matricial, there is a matricial subalgebra $R_0 \subseteq R/J$ that contains $\pi(\{1\} \cup X)$.  Then $R_1 := \pi^{-1}(R_0)$ is a unital subalgebra of $R$ containing $X$, and it suffices to show that $X$ is contained in a matricial subalgebra of $R_1$.  Moreover, $J$ is an ideal of $R_1$ and $R_1/J = R_0$ is matricial.  Thus, we may further reduce to the case where $R/J$ is matricial.

By Lemma \ref{R=T+J}, $R$ has a subalgebra $T$ such that $R = T \oplus J$ as vector spaces.  In particular, $T \cong R/J$ is matricial.  If $(t_1,\dots,t_n)$ is a $K$-basis for $T$, then $(\pi(t_1),\dots,\pi(t_n))$, where $\pi : R \rightarrow R/J$ is the quotient map, is a $K$-basis for $R/J$.

Now choose a finite subset $Y \subseteq J$ such that $X \subseteq T+Y$.  By Proposition \ref{q.cent.T}, there is an idempotent $e \in J$ such that $Y \subseteq eJe$ and $e$ centralizes $T$.  Since $1-e$ centralizes $T$, the set $S := (1-e)T$ is a subalgebra of $R$.  Moreover, $\pi$ sends $(1-e) t_i \mapsto \pi(t_i)$ for all $i$, hence  the $(1-e) t_i$ are $K$-linearly independent and $S \cong R/J \cong T$.  Further, the subspace $S+eJe$ is a subalgebra of $R$.
  
Since $Y \subseteq eJe$ and $eJe$ is locally matricial by Lemma \ref{lm>corner}, there is a matricial subalgebra $U \subseteq eJe$ which contains $Y \cup eT$, and $S+U$ is a matricial subalgebra of $R$.  Moreover, $T = (1-e)T + eT \subseteq S+U$, and therefore $X \subseteq T+Y \subseteq S+U$.
\end{proof}

Recall that \emph{ultramatricial} $K$-algebras are the direct limits of countable sequences of matricial algebras; these are just the locally matricial algebras with countable dimension.  Thus, Theorem \ref{lmxlm>lm} shows that the class of ultramatricial $K$-algebras is closed under extensions.

%%%%%%%%%%%%%%%%%%
\sectionnew{Locally f.d.ss by locally f.d.ss algebras}

We next prove both parts of Theorem B.  Theorem B(1), in parallel with Theorem \ref{lmxlm>lm}, will reduce to the case of an extension of a locally f.d.ss algebra by a f.d division algebra.  Again, several intermediate results are required.  After proving Theorem B, we show that two specialized types of extensions of locally f.d.ss algebras over arbitrary fields are locally f.d.ss.

\begin{lemma}  \label{lfdss>corner}
Suppose $R$ is a locally f.d.ss algebra, $J$ an ideal of $R$, and $e \in R$ an idempotent.  Then $eRe$ and $eJe$ are locally f.d.ss algebras.
\end{lemma}

\begin{proof}
As in Lemma \ref{lm>corner}, we only need to deal with $eJe$.

Let $X$ be a finite subset of $eJe$, and $f \in eJe$ an idempotent such that $X \subseteq f(eJe)f = fJf$.
By hypothesis, there is a f.d.ss subalgebra $T$ of $R$ that contains $X \cup \{f\}$, and it is enough to show that $fTf$ is f.d.ss.

Now $T = T_1 \oplus\cdots\oplus T_m$ where $T_i$ is an ideal of $T$ and $T_i \cong M_{n_i}(D_i)$ as algebras for some $n_i \in \NN$ and some f.d division algebras $D_i$.  Write $f = f_1+\cdots+f_m$ for some idempotents $f_i \in T_i$.  As a right $T_i$-module, $f_iT_i$ is a finite direct sum of some $r_i$ copies of the unique simple right $T_i$-module $S_i$, and so $f_iTf_i \cong M_{r_i}(\End_{T_i}(S_i)) \cong M_{r_i}(D_i)$.  Thus $fTf \cong M_{r_1}(D_1) \oplus\cdots\oplus M_{r_m}(D_m)$, proving that $fTf$ is f.d.ss, as required.
\end{proof}

\begin{lemma}  \label{lfdss.equiv0}
Let $R$ be an algebra, $J$ a locally f.d.ss ideal of $R$, and $T$ a f.d.ss subalgebra of $R$ such that $R = T \oplus J$ as vector spaces.  Then $R$ is locally f.d.ss if and only if
\begin{itemize}
\item[(*)] For any finite subset $Y$ of $J$, there is an idempotent $e \in J$ such that $Y \subseteq eJe$ and $e$ centralizes $T$.
\end{itemize}
\end{lemma}

\begin{proof}
$(\Longrightarrow)$: Since $R$ is locally f.d.ss, there is a f.d.ss subalgebra $U$ of $R$ containing $T\cup Y$. Since $U\cap J $ is an ideal of $U$, we must have $U \cap J = eU$ for some idempotent $e \in U\cap J$ which is central in $U$.  Then $Y \subseteq U\cap J = eUe \subseteq eJe$ and $e$ centralizes $T$.

$(\Longleftarrow)$:
Given a finite subset $X$ of $R$, there is a finite set $Y \subseteq J$ such that $X \subseteq T+Y$.  By (*), there is an idempotent $e \in J$ such that $Y \subseteq eJe$ and $e$ centralizes $T$.

Since $eT \cup Y \subseteq eJe$ and $eJe$ is locally f.d.ss (Lemma \ref{lfdss>corner}), there is a f.d.ss subalgebra $U$ of $eJe$ containing $\{e\}\cup eT\cup Y$.  Since $e$ centralizes $T$, the set $(1-e)T := \{ t-et \mid t \in T \}$ is a subalgebra of $R$.  The quotient map $R \rightarrow R/J$ sends both $T$ and $(1-e)T$ isomorphically onto $R/J$, so $(1-e)T$ is f.d.ss.  Finally, the subspace $V := (1-e)T+U$ is a f.d.ss subalgebra of $R$ containing $T+Y$, and $X \subseteq V$.
\end{proof} 

\begin{corollary}  \label{lfdss.equiv}
Let $R$ be a unital algebra, $J$ a locally f.d.ss ideal of $R$, and $T$ a f.d.ss subalgebra of $R$ with identity $p:=1_T$ such that $R = T \oplus J$ as vector spaces.  Then the following are equivalent:

{\rm(a)} $R$ is locally f.d.ss.

{\rm(b)} $pRp$ is locally f.d.ss.

{\rm(c)} For any idempotent $f \in pJp$, there is an idempotent $e \in pJp$ such that $f \in eJe$ and $e$ centralizes $T$.
\end{corollary}

\begin{proof}
Set $1 := 1_R$, and note that $1-p \in J$.

(a)$\Longrightarrow$(b) by Lemma \ref{lfdss>corner}.

(b)$\Longrightarrow$(c): Note that $pJp$ is an ideal of $pRp$ and $pRp = T \oplus pJp$ as vector spaces.  Moreover, $pJp$ is locally f.d.ss by Lemma \ref{lfdss>corner}.  Thus, (c) follows from Lemma \ref{lfdss.equiv0}.

(c)$\Longrightarrow$(a): By Lemma \ref{lfdss.equiv0}, it suffices to show that if $Y$ is any finite subset of $J$, there is an idempotent $e \in J$ such that $Y \subseteq eJe$ and $e$ centralizes $T$.  Choose an idempotent $g \in J$ such that $Y \cup \{1-p\} \subseteq gJg$.  Then $g$ commutes with $p$, so $pg$ is an idempotent in $pJp$, and $g = pg+1-p$.  By (c), there is an idempotent $h \in pJp$ such that $pg \in hJh$ and $h$ centralizes $T$.  Then $e := h+1-p$ is an idempotent in $J$ which centralizes $T$.  Since $pg \in hJh$, it follows that $g \in eJe$, and thus $Y \subseteq gJg \subseteq eJe$.
\end{proof}

\begin{remark}  \label{Kperf}
We require some classical facts about algebras over perfect fields. 

Recall that a f.d $K$-algebra $A$ is \emph{separable} (\emph{over $K$}) provided $A \otimes_K F$ is semisimple for all fields $F \supseteq K$.  If $E \supseteq K$ is a finite degree field extension, then $E$ is separable as a $K$-algebra if and only if it is separable as an extension field of $K$ (e.g., \cite[Prop 7.4]{CR}).

If $K$ is perfect, then every algebraic field extension of $K$ is separable (e.g., \cite[Corollary V.6.12]{Lan}), and hence every f.d.ss $K$-algebra is separable (e.g., \cite[Cor 7.6]{CR}).  Consequently, the class of f.d.ss $K$-algebras is closed under tensor products over $K$ (e.g., \cite[Cor 7.8(i)]{CR}).  It follows directly that the class of locally f.d.ss $K$-algebras is closed under tensor products over $K$.

The latter facts of course fail when $K$ is not perfect.  For instance, if $E \supsetneq K$ is a purely inseparable field extension of finite degree, then $E$ is a f.d.ss $K$-algebra but $E \otimes_K E$ is not semisimple.
\end{remark}

\begin{lemma}  \label{proj.over.tensor}
Assume that $K$ is perfect. Let $T \subseteq R$ be unital $K$-algebras such that $1_T = 1_R$ and $T$ is a f.d division algebra.  View the bimodule $_TR_R$ as a right module over the algebra $\Rtil := \top \otimes_K R$.  Then $R$ is a projective right $\Rtil$-module.
\end{lemma}

\begin{proof}
For emphasis, write the right $\Rtil$-module multiplication on $R$ with a dot, so that $r\cdot (t\otimes r') = trr'$ for $r,r' \in R$ and $t\in T$.

Set $\Ttil := \top \otimes_K T \subseteq \Rtil$ and view $T = {}_TT_T$ as a right $\Ttil$-module.  Since $\Ttil$ is f.d.ss (Remark \ref{Kperf}), $T_{\Ttil}$ is projective.  

Let $(b_\beta)_{\beta \in B}$ be a basis for $_TR$.  As left $\Ttil$-modules,
$$
\Rtil = \top \otimes_K \bigoplus_{\beta\in B} Tb_\beta = \bigoplus_{\beta\in B} (\top \otimes_K Tb_\beta) = \bigoplus_{\beta\in B} \Ttil(1\otimes b_\beta),
$$
whence $T \otimes_{\Ttil} \Rtil = \bigoplus_{\beta\in B} \bigl( T \otimes_{\Ttil} \Ttil(1\otimes b_\beta) \bigr) = \bigoplus_{\beta\in B} \bigl( T \otimes (1\otimes b_\beta) \bigr)$.

Let $\theta : T \otimes_{\Ttil} \Rtil \rightarrow R$ be the composition of the following right $\Rtil$-module homomorphisms:
$$
T \otimes_{\Ttil} \Rtil \xrightarrow{\;\text{incl}\,\otimes\,\text{id}\;} R \otimes_{\Ttil} \Rtil \xrightarrow{\;\text{canon}\;} R \otimes_{\Rtil} \Rtil \xrightarrow{\;\text{mult}\;} R,
$$
so that $\theta(t\otimes (s\otimes r)) = t\cdot(s\otimes r) = str$ for $t,s\in T$ and $r\in R$.
Since $\theta(1_T\otimes(1_T\otimes r)) = r$ for $r \in R$, we see that $\theta$ is surjective.

Given $x \in \ker \theta$, write $x = \sum_{\beta\in B} t_\beta \otimes (1\otimes b_\beta)$ for suitable $t_\beta \in T$, and observe that $\sum_{\beta\in B} t_\beta b_\beta = \theta(x) = 0$.
Then all $t_\beta = 0$, whence $x=0$.  Thus $T \otimes_{\Ttil} \Rtil \cong R$ as right $\Rtil$-modules via $\theta$.

The projectivity of $T_{\Ttil}$ therefore implies that $R_{\Rtil}$ is projective.
\end{proof}

\begin{theorem}  \label{perf.lfdss.ext}
Assume that $K$ is perfect.  Let $R$ be a $K$-algebra and $J$ an ideal of $R$.  If $J$ and $R/J$ are both locally f.d.ss $K$-algebras, then $R$ is locally f.d.ss.
\end{theorem}

\begin{proof}
{\bf Case 1:} $R$ is unital, $T$ is a f.d division subalgebra of $R$ with $1_T=1_R$, and $R = T\oplus J$ as vector spaces.

By Lemma \ref{lfdss.equiv0}, it suffices to show that for any finite subset $Y$ of $J$, there is an idempotent $e \in J$ such that $Y \subseteq eJe$ and $e$ centralizes $T$.
We first find an idempotent $f \in J$ such that $Y \subseteq fJ$ and $f$ centralizes $T$.

Let $\Rtil := \top \otimes_K R$ as in Lemma \ref{proj.over.tensor}, so that $R$ becomes a projective right $\Rtil$-module.  Note that $\top \otimes_K J$ is an ideal of $\Rtil$ with $\Rtil/(\top \otimes_K J) \cong \top \otimes_K T$.  By Remark \ref{Kperf}, $\top \otimes_K J$ and $\top \otimes_K T$ are locally f.d.ss $K$-algebras, hence regular.  Thus, $\Rtil$ is a regular ring.

Now $V := \sum_{y\in Y} TyR$ is a finitely generated right ideal of $R$, hence also a finitely generated right $\Rtil$-submodule of $R$.  Since $\Rtil$ is regular and $R_{\Rtil}$ is projective, $V$ is an $\Rtil$-module direct summand  of $R$ \cite[Theorem 1.11]{vnrr}, so $R = V \oplus W$ for some $\Rtil$-submodule $W$.  As $V$ and $W$ are also right ideals of $R$, there is an idempotent $f \in R$ such that $fR = V$ and $(1-f)R = W$.  Now $f \in J$ and $Y \subseteq V = fJ$ because $V \subseteq J$, while $T fR \subseteq fR$ and $T (1-f)R \subseteq (1-f)R$ because $V$ and $W$ are $\Rtil$ submodules of $R$.  The latter inclusions imply $(1-f)Tf = fT(1-f) = 0$, and therefore $f$ commutes with all elements of $T$.

Symmetrically, but with $Y \cup \{f\}$ in place of $Y$, there exists an idempotent $g \in J$ such that $Y \cup \{f\} \subseteq Jg$ and $g$ centralizes $T$.  Then $e := f+g-gf$ is an idempotent in $J$ such that $ef=f$ and $ge=g$, whence $Y \subseteq eJe$.  Since $f$ and $g$ centralize $T$, so does $e$.

{\bf Case 2:} $R$ is unital, $T$ is a f.d simple subalgebra of $R$ with $1_T=1_R$, and $R = T\oplus J$ as vector spaces.

Under the current assumptions, there is a set of $n\times n$ matrix units $e_{ij} \in T$, for some $n\in\NN$, such that $e_{11}+\cdots+e_{nn}=1$ and $T_1 := e_{11}Te_{11}$ is a f.d division algebra.  Note that $J_1 := e_{11}Je_{11}$  is an ideal of $R_1 := e_{11}Re_{11}$, that $J_1$ is locally f.d.ss (Lemma \ref{lfdss>corner}), that $1_{T_1} = 1_{R_1}$, and that $R_1 = T_1 \oplus J_1$ as vector spaces.  Thus, $R_1$ is locally f.d.ss by Case 1.

Given a finite subset $X$ in $R$, set $X_1 := \{ e_{1i}xe_{j1} \mid i,j=1,\dots,n,\ x \in X \} \subseteq R_1$.  Since $R_1$ is locally f.d.ss, there is a f.d.ss subalgebra $S_1$ of $R_1$ that contains $X_1$.  If $q_1 := 1_{S_1}$, then $K(e_{11}-q_1) + S_1$ is a unital f.d.ss subalgebra of $R_1$ containing $X_1$, and so we may assume that $q_1 = e_{11}$.   Now $S := \sum_{i,j=1}^n e_{i1}S_1e_{1j}$ is a unital subalgebra of $R$, isomorphic to $M_n(S_1)$ and hence f.d.ss.  Also,
$$
x = \sum_{i,j=1}^n e_{ii}xe_{jj} = \sum_{i,j=1}^n e_{i1}(e_{1i}xe_{j1})e_{1j} \in \sum_{i,j=1}^n e_{i1}X_1e_{1j} \subseteq S \qquad \forall\; x \in X.
$$

{\bf Case 3:} $R$ is unital and $R/J$ is f.d, simple, and unital.

In this case, Lemma \ref{R=T+J} implies that $R$ has a subalgebra $T$ such that $R = T\oplus J$ as vector spaces.  Hence, $T \cong R/J$ is a f.d, simple, unital algebra.
If $p := 1_T$, then $pJp$ is a locally f.d.ss ideal of $pRp$ (Lemma \ref{lfdss>corner}), $pRp/pJp \cong R/J$, and $pRp = T \oplus pJp$ as vector spaces.  Now $pRp$ is locally f.d.ss by Case 2, and hence Corollary \ref{lfdss.equiv} implies that $R$ is locally f.d.ss.

{\bf Case 4:} $R/J$ is f.d, simple, and unital.

Again, Lemma \ref{R=T+J} implies that $R$ has a subalgebra $T$ such that $R = T\oplus J$ as vector spaces, and $T \cong R/J$ via the quotient map $R \rightarrow R/J$.

Let $X$ be a finite subset of $R$.  Since $T$ is f.d, there is an idempotent $q \in R$ such that $T \cup X \subseteq qRq$.  Then $qJq$ is a locally f.d.ss ideal of $qRq$ (Lemma \ref{lfdss>corner}) and $qRq/qJq \cong R/J$.  By Case 3, $qRq$ is locally f.d.ss, and thus $qRq$ has a f.d.ss subalgebra $S$ that contains $X$.

{\bf Case 5:} $R/J$ is f.d.ss.

In this case, $R/J = A_1 \oplus\cdots\oplus A_n$ where the $A_j$ are ideals of $R/J$ which are f.d, simple, and unital as algebras.  Then $R$ has ideals $R_0 = J \subset R_1 \subset\cdots\subset R_n = R$ such that $R_j/R_{j-1} \cong A_j$ as algebras for $j=1,\dots,n$.  By Case 4, the local f.d.ss property passes from $R_{j-1}$ to $R_j$ for $j=1,\dots,n$, and thus $R$ is locally f.d.ss.

{\bf Case 6:} General case.

Let $X$ be an arbitrary finite subset of $R$, and let $\pi : R \rightarrow R/J$ be the quotient map.
Since $R/J$ is locally f.d.ss, there is a f.d.ss subalgebra $A$ of $R/J$ that contains $\pi(X)$.  Then $R' := \pi^{-1}(A)$ is a subalgebra of $R$ such that $X,J \subseteq R'$ and $R'/J = A$.  By Case 5, $R'$ is locally f.d.ss, and therefore $X$ is contained in a f.d.ss subalgebra of $R'$.
\end{proof}

To address Theorem B(2), we work within algebras of the following type.  For any field $L$, let $\rcfm(L)$ denote the $L$-algebra of $\omega\times\omega$ row- and column-finite matrices over $L$, and let $\fm(L)$ be the ideal of $\rcfm(L)$ consisting of matrices with at most finitely many nonzero entries.  Write $I_\omega$ for the identity element of $\rcfm(L)$. 

\begin{exampleproclaim}  \label{ex.charp}
If $K$ is not perfect, there exists a $K$-algebra $R$ with an ideal $J$ such that $J$ and $R/J$ are locally f.d.sss but $R$ is not.
\end{exampleproclaim}

\begin{proof}
Let $p := \chr(K)$ and $\alpha \in K\setminus K^p$.  Choose a field extension $L = K(\tau) \supset K$ with $\tau^p = \alpha$.  Since $X^p-\alpha \in K[X]$ is irreducible (e.g., \cite[Exer.~V.15]{Lan}), $\langle X^p-\alpha \rangle$ is a maximal ideal of $K[X]$ and $L \cong K[X]/\langle X^p-\alpha \rangle$.  We construct $R$ as a $K$-subalgebra of $\rcfm(L)$.  To simplify calculations, identify $\rcfm(L)$ with $\rcfm(M_p(L))$ in the natural way.

Choose a nonzero matrix $\nu \in M_p(L)$ such that $\nu^p = 0$, set 
$$
N := \begin{bmatrix}
\nu&\nu&0&0&0\\  0&\nu&\nu&0&0&\cdots\\
0&0&\nu&\nu&0\\
&\vdots&&\ddots&\ddots&\ddots  \end{bmatrix} \in \rcfm(L),
$$
and note that $N^p=0$.  Since $\tau I_\omega$ commutes with $N$, the matrix $t := \tau I_\omega + N$ satisfies $t^p = \alpha I_\omega$.  Consequently, the $K$-subspace $T := K I_\omega + K t +\cdots+ K t^{p-1}$ is a $K$-subalgebra of $\rcfm(L)$ isomorphic to $L$.

Now set $J := \fm(L)$ and $R := T+J$.  Then $R$ is a $K$-subalgebra of $\rcfm(L)$ and $J$ is an ideal of $R$.  Since $J$ is locally matricial as an $L$-algebra, it is locally f.d.ss as a $K$-algebra.  Also, $R/J \cong T \cong L$ as $K$-algebras, so $R$ is an extension of one locally f.d.ss $K$-algebra by another.

We claim that $R$ is not a locally f.d.ss $K$-algebra.  By Lemma \ref{lfdss.equiv0}, it suffices to exhibit a finite subset $Y \subseteq J$ for which there is no idempotent $e \in J$ such that $Y \subseteq eJe$ and $e$ centralizes $T$.  The latter condition amounts to $et=te$, which reduces to $eN=Ne$ because $ \tau I_\omega$ is central in $\rcfm(L)$.  We shall take $Y := \{y\}$ where 
$$
y := \begin{bmatrix}
I_p&0&0\\  0&0&0&\cdots\\  0&0&0\\  
&\vdots&&\ddots&  \end{bmatrix}.
$$

Suppose there is an idempotent $e \in J$ such that $y \in eJe$ and $eN=Ne$.  Choose an even positive integer $n$ large enough so that the nonzero entries of $e$ lie in the upper left $(n+1)\times(n+1)$ block.  Then
$$
e = \begin{bmatrix}
\begin{array}{c|ccc|cc}
I_p&0&\cdots&0&0&\cdots\\
\hline
0&b_{11}&\cdots&b_{1n}&0&\cdots\\  \vdots&\vdots&&\vdots&\vdots\\  0&b_{n1}&\cdots&b_{nn}&0&\cdots\\
\hline
0&0&\cdots&0&0\\  \vdots&\vdots&&\vdots&&\ddots 
\end{array}
\end{bmatrix}
\qquad \text{for some}\ b_{ij} \in M_p(L).
$$

Comparing the first rows of $eN$ and $Ne$, we obtain
\begin{equation}  \label{prow1}
\nu=\nu b_{11} \qquad\text{and}\qquad 0=\nu b_{1j}\ \ (j=2,\dots,n).
\end{equation}
Comparing rows $2$ through $n$, we obtain
\begin{equation}  \label{prowmid}
\begin{gathered}
b_{i1}\nu = \nu(b_{i1}+b_{i+1,1}) \qquad\text{and}\qquad b_{in}\nu=0\quad (i=1,\dots,n-1)  \\
(b_{i,j-1}+b_{ij})\nu = \nu(b_{ij}+b_{i+1,j})\quad (i=1,\dots,n-1;\; j=2,\dots,n).
\end{gathered}
\end{equation}
Finally, comparing $(n+1)$st rows, we obtain
\begin{equation}  \label{prowend}
\begin{gathered}
b_{n1}\nu = \nu b_{n1} \qquad\text{and}\qquad b_{nn}\nu = 0 \\
(b_{n,j-1}+b_{nj})\nu = \nu b_{nj}\quad (j=2,\dots,n).
\end{gathered}
\end{equation}

Equations \eqref{prowmid}, \eqref{prowend} imply
\begin{align*}
0 &= b_{in}\nu = - b_{i1}\nu + \sum_{j=2}^n (-1)^j (b_{i,j-1}+b_{ij})\nu = \sum_{j=1}^n (-1)^j \nu(b_{ij}+b_{i+1,j})  \quad (i=1,\dots,n-1)  \\
0 &= b_{nn}\nu = - b_{n1}\nu + \sum_{j=2}^n (-1)^j (b_{n,j-1}+b_{nj})\nu = \sum_{j=1}^n (-1)^j \nu b_{nj} \,.
\end{align*}
Consequently,
\begin{equation}  \label{prow1add}
 - \sum_{j=1}^n (-1)^j \nu b_{1j} = \sum_{i=1}^{n-1} (-1)^i \sum_{j=1}^n (-1)^j \nu(b_{ij}+b_{i+1,j}) + \sum_{j=1}^n (-1)^j \nu b_{nj} =0.
 \end{equation}

In view of \eqref{prow1}, however, \eqref{prow1add} implies $\nu=0$, which is false.

Because of this contradiction, there is no idempotent $e \in J$ such that $y \in eJe$ and $e$ centralizes $T$.  Therefore $R$ is not locally f.d.ss.  In fact, it follows from the proof of Lemma \ref{lfdss.equiv0} that $R$ has no f.d.ss $K$-subalgebra that contains $\{ I_\omega,t,y \}$.
\end{proof}

In spite of Example \ref{ex.charp}, there are restricted classes of locally f.d.ss algebras over non-perfect fields which are closed under extensions, the obvious one being the class of locally matricial algebras.  The following theorem of Chekanu leads to another class.

\begin{theorem}  \label{chekanu}
\cite{Che}
Every regular locally f.d algebra with bounded index of nilpotence is locally f.d.ss.
\end{theorem}

\begin{theorem}  \label{ext.via.chekanu}
Let $R$ be a $K$-algebra and $J$ an ideal of $R$.  If $J$ and $R/J$ are both locally f.d.ss with bounded index of nilpotence, then $R$ is locally f.d.ss with bounded index.
\end{theorem}

\begin{proof}
We have already noted that such an algebra $R$ is regular.  It is locally f.d by \cite[Proposition X.12.1]{Jac} and has bounded index by \cite[Proposition 7.7]{vnrr}.  Therefore Theorem \ref{chekanu} implies that $R$ is locally f.d.ss.
\end{proof}

There is a mixed case that we record as well:

\begin{theorem}  \label{lfdssxLM>lfdss}
Let $R$ be a $K$-algebra and $J$ an ideal of $R$.  If $J$ is locally f.d.ss and $R/J$ is locally matricial, then $R$ is locally f.d.ss.
\end{theorem}

\begin{proof}
We first make reductions like those in the proof of Theorem \ref{lmxlm>lm}.

Let $X$ be a finite subset of $R$, and $e \in R$ an idempotent such that $X \subseteq eRe$.  By Lemmas \ref{lfdss>corner} and \ref{lm>corner}, $eJe$ is locally f.d.ss and $eRe/eJe$ is locally matricial, and it suffices to show that $X$ is contained in a f.d.ss subalgebra of $eRe$.  Thus, we may assume that $R$ is unital.
Now let $\pi : R \rightarrow R/J$ be the quotient map.  Since $R/J$ is locally matricial, there is a matricial subalgebra $R_0 \subseteq R/J$ that contains $\pi(\{1\} \cup X)$.  Then $R_1 := \pi^{-1}(R_0)$ is a unital subalgebra of $R$ containing $X$, and it suffices to show that $X$ is contained in a matricial subalgebra of $R_1$.  Moreover, $J$ is an ideal of $R_1$ and $R_1/J = R_0$ is matricial.  Thus, we may reduce to the case where $R/J$ is matricial.

By Lemma \ref{R=T+J}, there is a matricial subalgebra $T$ in $R$ such that $R= T\oplus J$ as vector spaces.  Given any finite subset $Y \subseteq J$, Proposition \ref{q.cent.T} shows that there is an idempotent $e \in J$ such that $Y \subseteq eJe$ and $e$ centralizes $T$.  Consequently, Lemma \ref{lfdss.equiv0} implies that $R$ is locally f.d.ss.  Therefore $X$ is contained in a f.d.ss subalgebra of $R$.
\end{proof}

%%%%%%%%%%%%%%%%%%%%%%%%%%%%%
\sectionnew{Local unit-regularity}

There are several conditions that might be used to define a ring (resp., an algebra) $R$ being locally unit-regular.  The strongest would be to require that every finitely generated subring (resp., subalgebra) of $R$ is unit-regular.  This condition is very restrictive; for instance, it implies that $R$ has no nonzero nilpotent elements, i.e., $R$ is an abelian regular ring.  Consistent with the definitions in the Introduction, we define a ring (resp., algebra) $R$ to be \emph{locally unit-regular} provided every finite subset of $R$ is contained in a subring (resp., subalgebra) with identity which is unit-regular.  The ring and algebra cases may be described via a common condition on corners, as follows.

\begin{lemma}  \label{loc.ur}
Let $R$ be a ring or an algebra.  Then $R$ is locally unit-regular as a ring (resp., as an algebra) if and only if $R$ is regular and $eRe$ is unit-regular for all idempotents $e \in R$.
\end{lemma}

\begin{proof}
$(\Longleftarrow)$ follows from the fact that regular rings are locally unital, as noted in Section \ref{sec.regular}.

$(\Longrightarrow)$: In either case, any element $a \in R$ is contained in a regular subring of $R$, whence $a$ is regular.

Now let $e \in R$ be idempotent and $x \in eRe$.  By assumption, $R$ has a subring (resp., subalgebra) $S$ with identity which is unit-regular and contains $\{e,x\}$.  Then $eSe$ is unit-regular \cite[Corollary 4.7]{vnrr} and $x = exe \in eSe$.  There is a unit $u$ of $eSe$ such that $xux=x$.  Since $u$ is also a unit in $eRe$, this proves that $eRe$ is unit-regular.
\end{proof}

\begin{corollary}  \label{lur>corner}
Let $R$ be a locally unit-regular ring (resp., algebra), $J$ an ideal of $R$, and $e \in R$ an idempotent.  Then $eRe$ and $eJe$ are locally unit-regular.
\end{corollary}

\begin{proof}
If $f$ is an idempotent in $eJe$, then $f(eJe)f = fJf = fRf$, which is unit-regular since $R$ is locally unit-regular.  This shows that $eJe$ is locally unit-regular.  The case of $eRe$ follows by taking $J=R$.
\end{proof}

For arbitrary $K$, the class of locally unit-regular $K$-algebras fails to be closed under extensions, due to \cite[Example 1]{MeMo}.  In the unital case, there is a lifting condition on units (invertible elements) that determines which such extensions are unit-regular \cite[Lemma 3.5]{Bac}.  This criterion carries over to the general situation as follows.

\begin{proposition}  \label{criterion.ext.lur}
Let $R$ be an algebra with an ideal $J$ such that $J$ and $R/J$ are both locally unit-regular.  Then $R$ is locally unit-regular if and only if
\begin{itemize}
\item[$(\dagger)$] For each idempotent $e \in R$, every unit in $eRe/eJe$ lifts to a unit in $eRe$.
\end{itemize}
\end{proposition}

\begin{proof}
Apply \cite[Lemma 3.5]{Bac} to all corners $eRe$, for idempotents $e \in R$.
\end{proof}

There is a further condition, weaker than those above, which has been used in the literature to define local unit-regularity, as in \cite[Definition 6]{AbRa}.  We label this condition as follows:  A ring $R$ is \emph{elementwise locally unit-regular} provided that for each $x \in R$, there exist an idempotent $e \in R$ and a unit $u$ of $eRe$ such that $x \in eRe$ and $xux=x$.  Abrams and Rangaswamy showed that a unital ring is elementwise locally unit-regular if and only if it is unit-regular \cite[Lemma 3(1)]{AbRa}, hence locally unit-regular.  For non-unital rings, however, elementwise local unit-regularity is weaker than local unit-regularity:

\begin{exampleproclaim}  \label{weak.loc.ureg}
There exist elementwise locally unit-regular $K$-algebras which are not locally unit-regular.
\end{exampleproclaim}

\begin{proof}
Choose a unital regular algebra $S$ which is not unit-regular, such as the endomorphism algebra of an infinite dimensional vector space.  Then take $R := FM(S)$, the algebra of $\omega\times\omega$ matrices over $S$ with at most finitely many nonzero entries.  Since $e_{11}Re_{11} \cong S$, the algebra $R$ is not locally unit-regular.

We may view $R= \bigcup_{n=1}^\infty M_n(S)$, where each $M_n(S)$ is identified with the subalgebra of matrices in $R$ whose nonzero entries lie in the upper left $n\times n$ block.  Given $x \in M_n(S)$, there exists $y \in M_n(S)$ such that $xyx=x$.  Then $e := I_{2n} \in M_{2n}(S) \subset R$ is idempotent, $x = \left[\begin{smallmatrix} x&0\\ 0&0 \end{smallmatrix}\right] \in eSe$, and $u := \left[\begin{smallmatrix} y&I_n\\ I_n&0 \end{smallmatrix}\right]$ is a unit in $eRe$ (with inverse $\left[\begin{smallmatrix} 0&I_n\\ I_n&-y \end{smallmatrix}\right]$) such that $xux = x$.  This verifies that $R$ is elementwise locally unit-regular.
\end{proof}

Since f.d.ss algebas are unit-regular, all locally f.d.ss algebras are locally unit-regular.
In the presence of regularity, local finite dimensionality suffices, as the following proposition shows.  Recall that regular, locally f.d algebras need not be locally f.d.ss, by examples of May \cite{May}.

\begin{proposition}  \label{reg.lfd>lur}  
If $R$ is a regular, locally f.d algebra, then $R$ is locally unit-regular.
\end{proposition}

\begin{proof}
Given an idempotent $e \in R$ and a finite subset $X$ of $eRe$, there exists a f.d subalgebra $S$ of $R$ containing $X \cup \{e\}$.  Then $eSe$ is a f.d subalgebra of $eRe$ containing $X$, showing that $eRe$ is locally f.d.  Thus, there is no loss of generality in assuming that $R$ is unital.

Given $x \in R$, choose $y \in R$ such that $xyx=x$.  By assumption, there is a f.d subalgebra $A$ of $R$ that contains $\{1,x,y\}$.  Then $xA = xyA \cong yxA$ and $A = xA \oplus (1-xy)A = yxA \oplus (1-yx)A$.  By Krull-Schmidt cancellation, $(1-xy)A \cong (1-yx)A$, whence $\rann_R(x) = (1-yx)R \cong (1-xy)R \cong R/xR$.  This implies that $R$ is unit-regular (e.g., \cite[proof of Theorem 4.1]{vnrr}).
\end{proof}

Although extensions of one locally f.d.ss algebra by another need not be locally f.d.ss, they are at least locally unit-regular, as follows.

\begin{theorem}  \label{lfdssxlur>lur}
Let $R$ be a $K$-algebra and $J$ an ideal of $R$.  If $J$ is locally unit-regular and $R/J$ is locally f.d.ss, then $R$ is locally unit-regular.
\end{theorem}

\begin{proof}
Since $R$ is regular, we just need to show that $eRe$ is unit-regular for any idempotent $e \in R$.  Observe that $eJe$ is a locally unit-regular ideal of $eRe$ (Corollary \ref{lur>corner}) and that $eRe/eJe \cong (e+J)(R/J)(e+J)$ is locally f.d.ss (Lemma \ref{lfdss>corner}).  Thus, there is no loss of generality in assuming that $R$ is unital and $e = 1_R$.

Consider an arbitrary element $x \in R$, and let $\pi : R \rightarrow R/J$ be the quotient map.  By assumption, there is a f.d.ss subalgebra $A$ in $R/J$ such that $\pi(x) \in A$.  After replacing $A$ by $K (1_{R/J}-1_A) + A$ if necessary, we may assume that $A$ is a unital subalgebra of $R/J$.  Now $R' := \pi^{-1}(A)$ is a unital subalgebra of $R$ containing $J$, and $R'/J \cong A$.  Since $x \in R'$, it suffices to prove that $x$ is unit-regular in $R'$.  Thus, we may also assume that $R/J$ is f.d.ss.

By Lemma \ref{R=T+J}, there is a subalgebra $T$ of $R$ such that $R = T\oplus J$ as vector spaces.  Since $T \cong R/J$ via $\pi$, it is f.d.ss, hence unital and unit-regular.  Any unit $u$ in $R/J$ lifts to a unit $v$ in $T$, and then also to a unit $v+1_R-1_T$ in $R$.  Consequently, \cite[Lemma 3.5]{Bac} implies that $R$ is unit-regular, and therefore $x$ is unit-regular in $R$.
\end{proof}

%%%%%%%%%%%%%%%%%%%%%% References %%%%%%%%%%%%%%%%%%%%%%%%%%%%%%%%%%%%%%%

\end{document}